\input amstex
\documentstyle {amsppt}
\magnification=1200
\vsize=9.5truein
\hsize=6.5truein
\nopagenumbers 
\nologo



\topmatter

\title
Irreducibly represented groups
\endtitle

\author
Bachir Bekka and Pierre de la Harpe
\endauthor

\abstract
A group is {\it irreducibly represented}
if it has a faithful irreducible unitary representation.
For countable groups, a criterion for irreducible representability
is given, which generalises a result 
obtained for finite groups by W.~Gasch\"utz in 1954.
In particular,  torsionfree groups and infinite conjugacy class groups
are irreducibly represented. 

We indicate some consequences of this for operator algebras.
In particular, we charaterise up to isomorphism the countable
subgroups $\Delta$ of the unitary group 
of a separable infinite dimensional Hilbert space $\Cal H$
of which the bicommutants $\Delta ''$
(in the sense of the theory of von Neumann algebras)
coincide with the algebra of all bounded linear operators on $\Cal H$.

\endabstract



\thanks
The authors are grateful to the
Swiss National Science Foundation 
for its support.
\endthanks

\address
Bachir Bekka, UFR Math\'ematiques, Universit\'e de Rennes~1, Campus Beaulieu,
F--35042 Rennes Cedex.
E-mail: bachir.bekka\@univ-rennes1.fr
\endaddress

\address
Pierre de la Harpe, Section de Math\'ematiques, 
Universit\'e de Gen\`eve, C.P. 64, 
CH-1211 Gen\`eve 4. 
E-mail: Pierre.delaHarpe\@math.unige.ch
\endaddress


\endtopmatter

\document

\head{\bf
1.~Gasch\"utz Theorem for infinite groups, and consequences
}\endhead

Define a group  to be {\it irreducibly represented}
if it has a faithful irreducible unitary representation
and {\it irreducibly underrepresented} 
\footnote{
On the day of writing, Google shows 
29 000 000 entries for {\it represented groups,}
2 390 000 for {\it underrepresented groups,}
641 000 for {\it \lq\lq represented groups\rq\rq ,}
670 000 for {\it \lq\lq underrepresented groups\rq\rq ,}
and zero entry for {\it \lq\lq irreducibly underrepresented groups\rq\rq .}
In some sense at least, what we have to say is new.
}
if not.
For example, a finite abelian group is irreducibly represented
if and only if it is cyclic
(because finite subgroups of multiplicative groups of fields,
in particular finite subgroups of $\bold C^*$, are cyclic).
It is a straightforward consequence of Schur's lemma 
that a group  of which the centre
contains a non--cyclic finite subgroup is irreducibly underrepresented.
For finite groups, there are also standard examples 
of groups without centre which are irreducibly underrepresented
(see  Note~F in \cite{Burns--11});
moreover, there exists a criterion due to Gasch\"utz
who states for finite dimensional representations 
over algebraically closed fields of characteristic zero
the equivalence of Properties (i), (iv), and (v)
in Theorem~2 below
(see \cite{Gasch--54},
as well as  \cite{Huppe--88, \S \ 42} and  \cite{P\'alfy--79}).
\par

The purpose of the present paper is to extend Gasch\"utz' result 
to infinite groups and unitary representations; 
for the particular case of finite groups, our arguments provide
a new proof of the main result of \cite{Gasch--54}
(at least for complex representations).
For a generalisation of Gach\"utz' result of a rather different kind, 
see \cite{Tushe--93}.

Since our arguments use measure theory, 
it is convenient to avoid the difficulties
connected with non--standard spaces,
so that we assume systematically that 
{\it the groups involved are countable}
(see also Example~VII in Subsection 5.a below).
Moreover and from now on, 
{\it we write \lq\lq representation\rq\rq \ for
\lq\lq unitary representation\rq\rq \,}
and, similarly,
\lq\lq character\rq\rq \ for
\lq\lq unitary character\rq\rq .

\bigskip\bigskip

To formulate our results, we need the  following preliminaries.
Let $\Gamma$ be a group.
\par

Let $N$ be a normal subgroup of $\Gamma$. 
A representation $\sigma$ of $N$ is said to be {\it $\Gamma$--faithful} if
$\bigcap_{\gamma \in \Gamma} \ker (\sigma^{\gamma}) = \{e\}$, 
where $e$ denotes the unit element of the group and
where $\sigma^{\gamma}$ denotes the representation 
$n \longmapsto \sigma(\gamma n \gamma^{-1}$), 
namely the conjugate of $\sigma$ by $\gamma$.
For example, if $\bold V$ denotes the normal subgroup of order $4$ 
in the symmetric group $\operatorname{Sym}(4)$ on four letters, 
any character of $\bold V$ distinct from the unit character
is $\operatorname{Sym}(4)$--faithful 
(even though  $\bold V$ does not have any faithful character).
\par

If $\{S_i\}_{i \in I}$ is a family of subsets of $\Gamma$,
we denote by $\langle \{S_i\}_{i \in I} \rangle$
the subgroup of $\Gamma$ generated by $\bigcup_{i \in I}S_i$.
Following \cite{Remak--30}, we define a {\it foot} of $\Gamma$
to be a  {\it minimal} normal subgroup of $\Gamma$,
namely a normal subgroup $M$ in $\Gamma$ such that $M \ne \{e\}$, 
and any normal subgroup of $\Gamma$ contained in $M$
is either $M$ or $\{e\}$.
We denote by $\Cal F_{\Gamma}$ the set of finite feet of $\Gamma$.
The {\it minisocle} of $\Gamma$ is the subgroup $MS(\Gamma)$ of $\Gamma$
generated by the union of its finite feet;
it is a characteristic subgroup of $\Gamma$.
Let $\Cal A_{\Gamma}$ denote the subset of $\Cal F_{\Gamma}$
of abelian groups, 
and let $\Cal H_{\Gamma}$ be the complement of $\Cal A_{\Gamma}$
in $\Cal F_{\Gamma}$.
We define $MA(\Gamma)$ and $MH(\Gamma)$
to be the subgroups of $\Gamma$ generated by
$\bigcup_{A \in \Cal A_{\Gamma}}A$
and
$\bigcup_{H \in \Cal H_{\Gamma}}H$
respectively;
both are characteristic subgroups of $\Gamma$ contained in $MS(\Gamma)$.
By the usual convention, $MS(\Gamma) = \{e\}$ if $\Cal F_{\Gamma}$ is empty,
and similarly for $MA(\Gamma)$ and $MH(\Gamma)$.

\proclaim{1.~Proposition} Let $\Gamma$ be a group, 
and let the notation be as above.
\par

(i) Each $A \in \Cal A_{\Gamma}$ is isomorphic to $(\bold F_p)^n$
for some prime $p$ and some positive integer $n$ (depending on $A$).
\par

(ii) There exists a subset $\{A_i\}_{i \in I}$ of $\Cal A_{\Gamma}$
such that $MA(\Gamma) = \bigoplus_{i \in I} A_i$.
In particular, the group $MA(\Gamma)$ is abelian.
\par

(iii) For each $H \in \Cal H_{\Gamma}$,
the feet $S_1,\hdots,S_k$ of $H$ are conjugate in $\Gamma$, and simple. 
Moreover $H = S_1 \oplus \cdots \oplus S_k$.
\par

(iv) We have $MH(\Gamma) = \bigoplus_{H \in \Cal H_{\Gamma}} H$.
\par

(v) We have $MS(\Gamma) = MA(\Gamma) \oplus MH(\Gamma)$.
\endproclaim

For some examples of minisocles, see Section~5.a.
Here is our first main result.

\proclaim{2.~Theorem} Let $\Gamma$ be a countable group.
Let $MA(\Gamma) = \bigoplus_{i \in I} A_i$ and
\newline
$MS(\Gamma) = MA(\Gamma) \oplus MH(\Gamma)$ be as above.
The following properties are equivalent:
\roster
\item"(i)" 
$\Gamma$ is irreducibly represented;
\item"(ii)" 
$MA(\Gamma)$ has a $\Gamma$--faithful character;
\item"(iii)" 
$MS(\Gamma)$ has a $\Gamma$--faithful irreducible representation;
\item"(iv)" 
for every finite subset $E$ of $I$, 
there exists an element
$x_E$ in $MA_E(\Gamma) \Doteq \bigoplus_{i \in E}A_i$  such that
the $\Gamma$--conjugacy class of $x_E$ generates $MA_E(\Gamma)$;
\item"(v)"
for every pair of finite subsets $E \subset I$ and $F \subset \Cal H_{\Gamma}$,
there exists an element
$z_{E,F}$ in $MS_{E,F}(\Gamma) \Doteq \left( \bigoplus_{i \in E}A_i \right)
\oplus \left( \bigoplus_{H \in F}H \right)$ such that
the $\Gamma$--conjugacy class of $z_{E,F}$ generates $MS_{E,F}(\Gamma)$.
\endroster
In particular, a countable group $\Gamma$ 
has a faithful irreducible representation
as soon as $MA(\Gamma) = \{e\}$, 
and {\it a fortiori} as soon as $MS(\Gamma) = \{e\}$.
\endproclaim

The next corollary is a straightforward consequence of Theorem~2. 
Recall that a group is {\it icc}
if it is not reduced to one element and if
all its conjugacy classes distinct from $\{e\}$ are infinite.

\proclaim{3.~Corollary}
For a countable group to be irreducibly represented, 
any of the three following conditions is sufficient:
(i) the group is torsionfree, (ii) the group is icc, 
(iii) the group has a faithful primitive action
on an infinite set.
\endproclaim

The case of icc groups is well--known, sometimes with a different proof.
Indeed, a group is icc if and only if its von Neumann algebra
is a factor of type $II_1$ (Lemma 5.3.4 of \cite{RO--IV});
it is then a standard fact that the reduced C$^*$--algebra 
of an icc group
has a faithful irreducible representation,
so that {\it a fortiori} the group itself 
has a faithful irreducible  representation
(see for example Proposition 21 of \cite{Harpe--07}).
\par

For a group $\Gamma$ 
which has a faithful primitive action on an infinite set $X$
(see  \cite{GelGl}),
observe that any normal subgroup of $\Gamma$ not reduced to $\{e\}$
is transitive on $X$ and therefore infinite, 
so that $MS(\Gamma) = \{e\}$.

\bigskip 

Theorem~2 does not state anything on the dimensions
of the representations which can occur in (i).
Before providing some information, let us recall that
a group is {\it virtually abelian}
if it has an abelian subgroup of finite index.

\proclaim{4.~Theorem} For a countable group $\Gamma$,
the two following properties are equivalent:
\roster
\item"(i)" $\Gamma$ has an infinite dimensional faithful 
   irreducible representation;
\item"(ii)" $\Gamma$ has the properties of Theorem~2 
   and is not virtually abelian.
\endroster
In other words, the following properties are equivalent:
\roster
\item"(iii)" $\Gamma$ has a faithful irreducible representation,
and all its faithful irreducible representations are finite dimensional;
\item"(iv)" $\Gamma$ has the properties of Theorem~2 
   and is virtually abelian.
\endroster
\endproclaim 

Let $M$ be a von Neumann algebra. 
We denote by $\Cal U (M)$ 
the unitary group $\{X \in M \mid X^*X=XX^*=1\}$ of $M$
and by $S''$ the double commutant of a subset $S$ of $M$.
Recall that $M$ is a {\it factor} if its centre is reduced to $\bold C$,
a {\it factor of type I} if there exists a Hilbert space $\Cal H$
such that $M = \Cal L (\Cal H)$,
and a {\it factor of type I$_{\infty}$} 
in case $\Cal H$ is infinite dimensional
(moreover, we assume here that Hilbert spaces are separable).
For factors of type I, we write
$\Cal U(\Cal H)$ instead of $\Cal U (M)$.

\proclaim{5.~Corollary}
Let $M = \Cal L (\Cal H)$ be a factor of type I$_{\infty}$.
For a countable group $\Gamma$,
the following two properties are equivalent:
\roster
\item"---" there is a subgroup $\Delta$ of $\Cal U(\Cal H)$ isomorphic to
$\Gamma$ such that $\Delta '' = M$;
\item"---" $\Gamma$ has the properties of Theorem~2 
and is not virtually abelian.
\endroster
\endproclaim

It would be interesting to have some information
of this kind for other factors.
In particular, we do not know any analogue 
of Theorem~4 
for any given finite dimension~$n \ge 2$,
nor of Corollary~5
for the finite dimensional factor $\Cal L (\bold C^n)$.
We do not know any solution to the {\it a priori} easier problem
to characterise the countable groups
which have at least 
one finite dimensional faithful irreducible representation.

\bigskip

The proof of Proposition 1 uses standard arguments
(compare with Section 4.3 of \cite{DixMo--96}). 
For the convenience of the reader, we give details in Section~2.
Theorem~2 is proved in Section~3. 
Theorem~4 and Corollary~5 are proved in Section~4.
We formulate a few remarks in Section~5:
on examples of socles and minisocles,
on the comparison between minisocles and periodic FC--kernels,
on a theorem of Gelfand and Raikov,
on tensor products of faithful representations,
and on countable groups with primitive maximal C$^*$--algebras.
The final Section~6 is devoted to
a generalisation of Theorem~2
concerning a countable group $\Gamma$
given together with a group of automorphisms $G$
which contains the group of inner automorphisms.

\bigskip

Understanding groups of a given class 
includes understanding their faithful actions of various kinds,
and the setting of linear (or unitary) actions is only one 
among several others.
For example, in the case of finite groups,
the questions of classifying 
multiply transitive actions and primitive actions
which are faithful have been central in group theory 
for more than hundred years;  
faithful primitive actions for infinite groups 
have been addressed in \cite{GelGl}.
Faithful amenable actions are the subject of \cite{GlaMo--07}.
Our initial motivation has been to ask 
some of the corresponding questions for linear actions.
\par

We are most grateful to Yair Glasner for explaining us his work \cite{GelGl}
and for his contribution to the setting out of the present work,
to Yehuda Shalom for a useful observation, 
and to Yves de Cornulier and John Wilson for their remarks 
on a preliminary version of this paper.

\bigskip
\head{\bf
2.~Proof of Proposition~1
}\endhead

We prepare the proof of Proposition~1 by recalling two lemmas.

\proclaim{6.~Lemma}
Let $\Gamma$ be a group. Let $M$ be a minimal normal subgroup of $\Gamma$
and $N$ a normal subgroup of $\Gamma$.
Then either $M \subset N$ or $\langle M,N \rangle = M \oplus N$.
\endproclaim

\demo{Proof} 
We can assume $N \ne \{e\}$.
Since $M \cap N$ is both in $M$ and normal in $\Gamma$,
either $M \cap N = M$, and $M \subset N$,
or $M \cap N = \{e\}$, and $\langle M,N \rangle = M \oplus N$.
\hfill $\square$
\enddemo

\proclaim{7.~Lemma}
Let $A$ be a group
and let $\left( S_i \right)_{i \in I}$ be a family of nonabelian simple groups;
set $S = A \oplus \left( \bigoplus_{i \in I} S_i \right)$.
Let $M$ be a minimal normal subgroup of $S$.
\par

Then either $M = S_{\ell}$ for some $\ell \in I$, or $M \subset A$.
\endproclaim

\demo{Proof}
Assume that $M \ne S_{\ell}$ for all $\ell \in I$.
Choose $i \in I$;
by Lemma~6 applied to $M$ and $N = S_i$, 
the groups $M$ and $S_i$ commute.
It follows that $M$ is a subgroup of the centraliser
of $\bigoplus_{i \in I} S_i$ in $S$,
namely a subgroup of $A$.
\hfill $\square$
\enddemo

\subhead 
Proof of Proposition~1 
\endsubhead 
(i) Let $A \in \Cal A_{\Gamma}$.
By the structure theory of finite abelian groups,
there exist a prime $p$ and an element $a \in A$ of order $p$.
Let $A^*$ denote the set of elements of order $p$ in $A$.
Then $A^* \cup \{e\}$ is a characteristic subgroup of $A$,
and therefore a normal subgroup of $\Gamma$.
By minimality of $A$, we have $A^* \cup \{e\} = A$,
so that $A$ is isomorphic to $(\bold F_p)^n$ 
for some $n \ge 1$, as claimed. 

\medskip

(ii) Let $\Cal L$ be 
the set of subsets $\left\{ A_{\ell} \right\}_{\ell \in L}$
of $\Cal A_{\Gamma}$ such that
$\langle \left\{ A_{\ell} \right\}_{\ell \in L} \rangle
= \bigoplus_{\ell \in L} A_{\ell}$;
we order $\Cal L$ by inclusion.
The crucial observation is that the ordered set $\Cal L$ is inductive,
so that we can choose a maximal element, 
say $\left\{ A_{i} \right\}_{i \in I}$.
Suppose that $\bigoplus_{i \in I}A_i$ is strictly contained in $MA(\Gamma)$;
we will arrive at a contradiction.
\par

Choose $B \in \Cal A_{\Gamma}$ 
such that $B$ is not contained in $\bigoplus_{i \in I}A_i$.
By Lemma~6 applied to $M = B$ and $N = \bigoplus_{i \in I}A_i$,
we have either $B \subset \bigoplus_{i \in I}A_i$,
which is ruled out by the choice of $B$,
or $\langle B, \left\{ A_i \right\}_{i \in I} \rangle = 
B \oplus \left( \bigoplus_{i \in I}A_i \right)$,
which is ruled out by the maximality of $I$.
This is the announced contradiction.

\medskip

(iii) Let $H \in \Cal H_{\Gamma}$.
Choose a minimal normal subgroup $S$ in $H$
(this is possible since $H$ is finite).
For each $x \in \Gamma$, the subgroup $x S x^{-1}$ is minimal normal in $H$.
Choose a set $S_1,\hdots,S_k$ of such conjugates of $S$ in $\Gamma$
which is such that 
$\langle S_1, \hdots, S_k \rangle = S_1 \oplus \cdots \oplus S_k$
and which is maximal for this property.
Set $N = \langle S_1, \hdots, S_k \rangle$;
it is a normal subgroup of $H$.
\par

We claim that $x S s^{-1} \subset N$ for each $x \in \Gamma$,
so that $N$ is normal in $\Gamma$.
Indeed, by Lemma~6 applied to $M = x S x^{-1}$ and $N$ in $H$,
either $\langle x S x^{-1}, S_1, \hdots, S_k \rangle
= x S x^{-1} \oplus S_1 \oplus \cdots \oplus S_k$,
but this is ruled out by the maximality of the set $\{S_1, \hdots, S_k\}$,
or $x S x^{-1} \subset N$, 
and this establishes the claim.
\par

Since $N$ is normal in $\Gamma$ and $N \subset H$,
we have $N = H$ by minimality of $H$.
Observe that, for each $i \in \{1,\hdots,k\}$, 
any normal subgroup of $S_i$ is normal in $H$;
it follows that $S_i$ is a simple group.
Finally, the set $\{S_1, \hdots, S_k\}$ coincides with
the set of all minimal normal subgroups of $H$ by Lemma~7.

\medskip

(iv) The same argument as for (ii) shows that 
there exists a subset $\left\{ H_k \right\}_{k \in K}$ of $\Cal H_{\Gamma}$
such that $\bigoplus_{k \in K} H_k = MH(\Gamma)$,
and Lemma~7 implies that 
$\left\{ H_k \right\}_{k \in K} = \Cal H_{\Gamma}$.

\medskip

(v) Again by the same argument as for (ii),
there exists a subset 
$\left\{ M_{\ell} \right\}_{\ell \in L}$ of $\Cal F_{\Gamma}$
such that $\bigoplus_{\ell \in L} M_{\ell} = MS(\Gamma)$,
and Lemma~7 implies that 
$\left\{ M_{\ell} \right\}_{\ell \in L}$ contains $\Cal H_{\Gamma}$.
\hfill $\square$

\bigskip
\head{\bf
3.~Proof of Theorem 2
}\endhead

We will prove successively that
$$
\matrix
(i) & \Longrightarrow & (ii) \hskip.1cm \& \hskip.1cm (iii) & 
   \text{(see Lemma 9),} \\
(iii) & \Longrightarrow & (i) &  \text{(Lemma 10),} \\
(ii) & \Longleftrightarrow & (iii) &  \text{(Lemma 13),} \\
(ii) & \Longleftrightarrow & (iv) &  \text{(Lemma 14).}
\endmatrix
$$
The equivalence $(iv) \Longleftrightarrow (v)$ is straightforward,
since nonabelian feet are direct products of simple groups.
Recall that we write  \lq\lq representation\rq\rq \ for
\lq\lq unitary representation\rq\rq .
\par

Given a representation $\pi$ of a countable  group $\Gamma$ 
in a Hilbert space~$\Cal H$, there exist
a standard Borel space $\Omega$, 
a bounded positive measure $\mu$ on $\Omega$,
a measurable field $\omega \longmapsto \pi_{\omega}$
   of irreducible  representations of $\Gamma$
in a measurable field $\omega \longmapsto \Cal H_{\omega}$
   of Hilbert spaces on $\Omega$,
and an isomorphism of $\Cal H$ 
   with $\int^{\oplus}_{\Omega} \Cal H_{\omega} d \mu (\omega)$
which implements a unitary equivalence
$$
\pi(\gamma) \, \approx \, \int^{\oplus}_{\Omega} \pi_{\omega}(\gamma) d\mu (\omega)
$$
for all $\gamma \in \Gamma$.
See \cite{Dix--69C$^*$, Sections 8.5 and 18.7.6}.
(Such decompositions in irreducible representations carry over to
continuous  representations of separable locally compact groups,
and more generally of separable C$^*$--algebras.
They are applications of the {\it reduction theory} for von Neumann algebras
\cite{Dix--69vN, Chapter~II}).
The following lemma is standard,
but we haven't found any appropriate reference.

\proclaim{8.~Lemma}
Let $\Gamma$ be a countable group.
Let $\Omega$ be a measure space with a  positive measure~$\mu$.
Let $\omega \longmapsto \pi_{\omega}$
be a measurable field of representations of $\Gamma$
in a measurable field of Hilbert spaces
$\omega \longmapsto \Cal H_{\omega}$
over $\Omega$
and let $\gamma \in \Gamma$.
\par
Then
$ \{\omega \in \Omega \mid \pi_{\omega} (\gamma) = I \}$
is a measurable subset of $\Omega$.
\endproclaim

\demo{Proof}
Let $(\xi^{(1)},\xi^{(2)}, ....)$ be
a fundamental sequence of measurable vector fields
(see \cite{Dix-69vN, Chapter~II, Number~1.3}).
For $i,j \geq 1$, consider the set
$$
\Omega_{i,j} \, = \, \{
\omega \in \Omega  \mid
\langle \pi_{\omega}(\gamma)\xi^{(i)}(\omega), \xi^{(j)}(\omega) \rangle =
\langle \xi^{(i)}(\omega), \xi^{(j)}(\omega) \rangle
\}.
$$
Observe that
$$
\{\omega \in \Omega \mid \pi_{\omega} (\gamma) = I \}  
\, = \, \bigcap_{i,j\geq 1}\Omega_{i,j}.
$$
Therefore it suffices to show that each set $\Omega_{i,j}$ is measurable.
\par
For fixed  $i,j \geq 1,$ the functions
$$
\omega \, \longmapsto \,
\langle \xi^{(i)}(\omega), \xi^{(j)}(\omega) \rangle
\qquad \text{and} \qquad
\omega \, \longmapsto \,
\langle \pi_{\omega}(\gamma)\xi^{(i)}(\omega), \xi^{(j)}(\omega) \rangle
$$
are measurable, by definition of a measurable vector field and of a measurable
field of representations. Hence $\Omega_{i,j}$ is measurable
and the  proof is complete.
\hfill $\square$
\enddemo

Let us now recall a general fact 
which can be seen as a weak form of Clifford theorem
for infinite dimensional representations.
(For a version of Clifford theorem concerning 
finite dimensional representations but possibly infinite groups,
see  Theorem 2.2 in \cite{Dixon--71}.)

\proclaim{9.~Lemma} Let $\Gamma$ be a countable group, 
$N$ a normal subgroup, 
$\pi$ an irreducible representation of~$\Gamma$ in a Hilbert space $\Cal H$, 
and $\sigma$ the restriction of $\pi$ to $N$.
Identify $\sigma$ to a direct integral of irreducible representations
$$
\sigma \, = \, \pi\vert_{N} \, = \, 
\int^{\oplus}_{\Omega} \sigma_{\omega} d\mu (\omega)
$$
as above.
\par

If  the representation $\pi$ is faithful,
then the representation $\sigma_{\omega}$ is $\Gamma$--faithful
for almost all $\omega \in \Omega$.
\endproclaim

\demo{Proof} If $N = \{e\}$, there is nothing to prove.
We assume from now on that $N$ is not reduced to one element.
\par

Denote by $\left\{C_j\right\}_{j \in J}$ 
the family of $\Gamma$--conjugacy classes in $N$ distinct from $\{e\}$.
For each $j \in J$, denote by $N_j$
the subgroup of $N$ generated by $C_j$;
observe that each $N_j$ is normal in $\Gamma$,
and that the family $\left\{ N_j \right\}_{j \in J}$
is countable (possibly finite) and nonempty. Set
$$
\Omega_j \, = \, \Big\{ 
\omega \in \Omega \mid
N_j \subset \ker \Big(
\bigoplus_{\gamma \in \Gamma} \hskip.1cm \sigma_{\omega}^{\gamma} 
\Big)
\Big\} 
\quad \text{and} \quad
\widetilde \Omega \, = \, \bigcup_{j \in J} \Omega_j \ .
$$
For $\omega \in \Omega$, 
observe that $\sigma_{\omega}$ is not $\Gamma$--faithful
if and only if the kernel of 
$\bigoplus_{\gamma \in \Gamma} \sigma_{\omega}^{\gamma}$
contains one of the $N_j$;
thus $\widetilde \Omega$ is the subset of $\Omega$
of the points $\omega$ such that $\sigma_{\omega}$ is not $\Gamma$--faithful.
Each $\Omega_j$ is measurable in $\Omega$
as a consequence of Lemma~8;
as $J$ is countable,  $\widetilde \Omega$ is also measurable.
\par

To end the proof, we assume that $\mu(\widetilde \Omega) > 0$
and we will arrive at a contradiction.
\par

As the family $J$ is countable, there exists $\ell \in J$ such that
$\mu(\Omega_{\ell}) > 0$.
Hence the unit representation $1_{N_{\ell}}$ of the group $N_{\ell}$
is strongly contained in the restriction of $\sigma$ to $N_{\ell}$,
so that the subspace of $\Cal H$ of $N_{\ell}$--invariant vectors
is not reduced to $\{0\}$.
Since $N_{\ell}$ is normal in $\Gamma$, 
this subspace is invariant by $\pi(\Gamma)$; 
by irreducibility of~$\pi$, 
this subspace is the whole of $\Cal H$. 
In other words,  the restriction of $\pi$ to $N_{\ell}$ 
is the unit representation. 
The last statement is a contradiction, since $\pi$ is faithful.
\hfill $\square$ 
\enddemo

The particular case of Lemma~9 
for which $N = MA(\Gamma)$ [respectively $N = MS(\Gamma)$] 
shows that (i) implies (ii) [respectively (iii)] in Theorem~2.
The implication $(iii) \Longrightarrow (i)$
follows from the next lemma applied to $N = MS(\Gamma)$ since,
by definition, there does not exist any finite foot $M$ of $\Gamma$
such that $M \cap MS(\Gamma) = \{e\}$.

\proclaim{10.~Lemma} Let $\Gamma$ be a countable group,
$N$ a normal subgroup,
$\sigma$ an irreducible  representation of $N$
in a Hilbert space $\Cal K$,
and $\pi = \operatorname{Ind}_N^{\Gamma}(\sigma)$
the corresponding induced representation.
Let $\pi = \int^{\oplus}_{\Omega} \pi_{\omega} d\mu(\omega)$
be a direct integral decomposition of $\pi$ 
into irreducible representations.
Assume that there does not exist any finite foot $M$ in $\Gamma$ 
such that $M \cap N = \{e\}$.
\par

If the representation $\sigma$ is $\Gamma$--faithful,
then the representation $\pi_{\omega}$ is faithful
for almost all $\omega$ in $\Omega$.
\endproclaim 

\demo{Proof}
In the model we choose for induced representations,
$\pi$ acts on the Hilbert space $\Cal H$
of mappings $f : \Gamma \longrightarrow \Cal K$ with the two following
properties:
$$
\aligned
&(1) \hskip.5cm f(\gamma n) \, = \, \sigma(n^{-1})f(\gamma)
   \hskip.5cm \text{for all} \hskip.2cm \gamma \in \Gamma
   \hskip.2cm \text{and} \hskip.2cm 
    n \in N, \hskip1.5cm
\\
&(2) \hskip.5cm \sum_{
              \Gamma/ N
                     }
   \Vert f(\gamma) \Vert^2 \, < \, \infty. \hskip1.5cm
\endaligned
$$
(The notation of (2) indicates a summation over
one representant $\gamma \in \Gamma$ of each class in $\Gamma/ N$.)
Then $\left(\pi(x)f\right)(\gamma) = f(x^{-1}\gamma)$
for all $x, \gamma \in \Gamma$.
\par

Denote this time by $\{C_j\}_{j \in J}$ 
the family of conjugacy classes of $\Gamma$ 
distinct from $\{e\}$.
For each $j \in J$, denote by $\Gamma_j$ the subgroup generated by $C_j$,
which is a normal subgroup of $\Gamma$ not reduced to $\{e\}$; 
set
$$
\Omega_j \, = \, \left\{ 
\omega \in \Omega \mid
\Gamma_j \subset \ker \left(\pi_{\omega} \right)
\right\} 
\quad \text{and} \quad
\widetilde \Omega \, = \, \bigcup_{j \in J} \Omega_j .
$$
As in the proof of Lemma~9, 
$\widetilde \Omega$ is the set of points $\omega$
such that $\pi_{\omega}$ is not faithful,
and it is measurable.
To end the proof, we assume that $\mu(\widetilde \Omega) > 0$,
so that there exists $\ell \in J$ for which $\mu(\Omega_{\ell}) > 0$,
and we will arrive at a contradiction.
\par

Continuing as in the proof of Lemma~9,
we observe that
there exists a nonzero vector $f : \Gamma \longrightarrow \Cal K$ in 
$\Cal H = \int^{\oplus}_{\Omega} \Cal H_{\omega} d\mu(\omega)$
which is supported in $\Omega_{\ell}$
(as a measurable section of the field of Hilbert spaces
$\omega \longmapsto \Cal H_{\omega}$
underlying the field of representations
$\omega \longmapsto \pi_{\omega}$),
and which is such that $\pi(x)f = f$ for all $x \in \Gamma_{\ell}$.
\par

Let $\gamma_0 \in \Gamma$ be such that $f(\gamma_0^{-1}) \ne 0$;
set $\xi = f(\gamma_0^{-1})$.
Using (1), we find 
$$
(3) \hskip.5cm 
\xi \, = \, 
f(\gamma_0^{-1}) \, = \,
f(x^{-1}\gamma_0^{-1}) \, = \,
f \left( \gamma_0^{-1} (\gamma_0 x^{-1} \gamma_0^{-1}) \right) \, = \,
\sigma(\gamma_0 x \gamma_0^{-1})\xi \, = \,
\sigma^{\gamma_0} (x) \xi
$$
for all $x \in \Gamma_{\ell} \cap N$.

\medskip

{\it Claim~1: $\Gamma_{\ell} \cap N = \{e\}$}.
Denote by $\Cal K^{\Gamma_{\ell} \cap N}$ the subspace of $\Cal K$
of vectors invariant by $\sigma^{\gamma_0} (\Gamma_{\ell} \cap N)$.
This is a $\sigma^{\gamma_0}(N)$--invariant
subspace of $\Cal K$, 
since $\Gamma_{\ell} \cap N$ is a normal subgroup of $N$.
Now $\Cal K^{\Gamma_{\ell} \cap N} \ne \{0\}$ by (3)
and $\Cal K^{\Gamma_{\ell} \cap N} = \Cal K$ 
because $\sigma^{\gamma_0}$ is irreducible.
Thus $\Gamma_{\ell} \cap N$ is inside the kernel
of the representation $\sigma^{\gamma_0}$ of $N$;
as $\Gamma_{\ell} \cap N$ is normal in $\Gamma$,
the group $\Gamma_{\ell} \cap N$ is also inside the kernel 
of the representation $\sigma^{\gamma}$ of $N$
for all $\gamma \in \Gamma$.
As $\sigma$ is $\Gamma$--faithful, $\Gamma_{\ell} \cap N = \{e\}$, as
claimed.

\medskip

{\it Claim~2: the subgroup $\Gamma_{\ell}$ of $\Gamma$ is finite.}
Consider the function
$$
\varphi : \Gamma \longrightarrow \bold R_+, \hskip.5cm
\gamma \longmapsto \Vert f(\gamma) \Vert .
$$
We have
$$
\aligned
&(4) \hskip.5cm 
\varphi(\gamma_0^{-1}) \, \ne \, 0,
\\
&(5) \hskip.5cm 
\varphi \quad 
   \text{is constant under right translations
        by elements of $N$},
\\
&(6) \hskip.5cm 
\sum_{\Gamma / N} \vert \varphi(\gamma) \vert^2 \, < \, \infty,
\\
&(7) \hskip.5cm 
\varphi \quad
   \text{is invariant under left translations by elements of $\Gamma_{\ell}$}.
\endaligned
$$
It follows from (4) to (7) that the image of $\Gamma_{\ell}\gamma_0$
in $\Gamma / N$ is finite.
The image of $\gamma_0^{-1} \Gamma_{\ell} \gamma_0 = \Gamma_{\ell}$
in $\Gamma / N$ is also finite,
so that the index of $N$ in $\Gamma_{\ell} N$ is finite.
Claim~2 follows since 
$\Gamma_{\ell} N$ is isomorphic to the direct sum
$\Gamma_{\ell} \oplus N$ by Claim~1.

\medskip

Any  subgroup $M$ of $\Gamma_{\ell}$ which is normal in $\Gamma$
and minimal for this property is a finite foot of $\Gamma$,
and $M \cap N = \{e\}$ by Claim~1.
This is in contradiction with one of the hypotheses of the lemma.
\hfill $\square$ 
\enddemo

The particular case $N = \{e\}$ is of independent interest.

\proclaim{11.~Proposition} Let $\Gamma$ be a countable infinite group
which does not contain any finite foot, and let
$
\lambda_{\Gamma} = \int^{\oplus}_{\Omega} \pi_{\omega} d\mu(\omega)
$
be a direct integral decomposition of $\lambda_{\Gamma}$ 
into irreducible representations.
Then $\pi_{\omega}$ is faithful for almost all $\omega \in \Omega$.
\endproclaim

Next, we show that $(ii) \Longleftrightarrow (iii)$ in Theorem~2. 
This will be a consequence of Lemma~13,
for the proof of which we will call upon the following lemma.
\par

For a Hilbert space $\Cal H$, we denote by $\Cal L (\Cal H)$
its algebra of bounded linear operators.

\proclaim{12.~Lemma}
Let $\Cal H_1,\Cal H_2$ be two Hilbert spaces.
Let $S_1\in \Cal L(\Cal H_1), S_2\in \Cal L(\Cal H_2)$
be such that
$S_1 \otimes S_2 \in \Cal L(\Cal H_1 \otimes \Cal H_2)$
is a non--zero multiple of the identity operator. 
\par

Then $S_1$ and $S_2$ are multiples of the identity.
\endproclaim

\demo{Proof}
Let $\lambda \in \bold C^*$ be such that 
$S_1 \otimes S_2= \lambda I$.
Let $\{\xi_i\}_{i\in I}$ 
be a Hilbert space basis of $\Cal H_1$. 
Since $S_2 \neq 0$, 
there exist $\eta_1,\eta_2\in\Cal H_2$ such that 
$$
\langle S_2(\eta_1),\eta_2\rangle\neq 0.
$$
For every $\xi\in \Cal H_1,$  we have
$$
\langle(S_1\otimes S_2)(\xi\otimes \eta_1), \xi_i\otimes \eta_2\rangle
=\langle S_1(\xi), \xi_i\rangle\langle S_2(\eta_1), \eta_2\rangle
$$
and hence
$$
\aligned
\langle S_1(\xi), \xi_i\rangle
& = \dfrac{1}{\langle S_2(\eta_1), \eta_2\rangle}
\langle(S_1\otimes S_2)(\xi \otimes \eta_1), \xi_i \otimes \eta_2\rangle
\\
& = \dfrac{\lambda}{\langle S_2(\eta_1), \eta_2\rangle}
\langle\xi\otimes \eta_1, \xi_i\otimes \eta_2\rangle
\\
& = \dfrac{\lambda\langle \eta_1, \eta_2\rangle}{\langle S_2(\eta_1), \eta_2\rangle}
\langle\xi, \xi_i\rangle\\
\endaligned
$$
for all $i\in I.$ It follows that
$$
\aligned
S_1(\xi)
& = \sum_{i\in I} \langle S_1(\xi), \xi_i\rangle \xi_i
\\
& = \dfrac{\lambda\langle \eta_1, \eta_2\rangle}{\langle S_2(\eta_1),
\eta_2\rangle}\sum_{i\in I}
\langle \xi, \xi_i\rangle \xi_i
\\
& = \dfrac{\lambda\langle \eta_1, \eta_2\rangle}{\langle S_2(\eta_1), \eta_2\rangle}
\xi
\endaligned
$$
for every $\xi\in H_1,$ showing that $S_1$ is a multiple of the identity.
A similar argument applies to  $S_2$. 
\hfill $\square$ 
\enddemo

\proclaim{13.~Lemma} 
Let $\Gamma$ be a group
and let $N$ be a normal subgroup of $\Gamma$.
Assume that $N = A \oplus S$, 
where $A$ is an abelian normal subgroup of $\Gamma$ and 
where $S$ is the direct sum of a family $\left( S_i \right)_{i \in I}$ of
finite simple nonabelian normal subgroups of $S$.
The following properties are equivalent:
\roster
\item"(i)"
$N$ has a $\Gamma$-faithful irreducible representation;
\item"(ii)"
$A$ has a $\Gamma$-faithful character.
\endroster
\endproclaim

\demo{Proof}
Assume first that there exists 
a $\Gamma$--faithful irreducible representation $\pi$ of $N$.
Since the factor $A$ of $N = A \oplus S$ is abelian,
and in particular a type~I group,
there exist a character $\chi$ of $A$
and an irreducible representation $\rho$ of $S$
such that $\pi = \chi \otimes \rho$
\cite{Dix--69C$^*$, Proposition 13.1.8}.
Since $\ker (\chi^{\gamma}) = \ker\left( (\pi^{\gamma}) \vert_A \right)$
for all $\gamma \in \Gamma$,
the character $\chi$ of $A$ is $\Gamma$--faithful.
\medskip

Assume now that there exists 
a $\Gamma$--faithful character $\chi$ of $A$.
We claim that there exists an irreducible representation $\rho$ of $S$
such that, for every $\gamma \in S$, $\gamma \ne e$,
the operator $\rho(\gamma)$ is not a multiple of the identity operator.
Lemma~12 will then imply that
the exterior tensor product $\chi \otimes \rho$
is a $\Gamma$--faithful representation of $N = A \oplus S$.
\par

For every $i \in I$, 
let $\rho_i$ be an irreducible representation  of $S_i$ 
distinct from the unit representation, in some Hilbert space~$\Cal H_i$.
Choose a unit vector $\eta_i \in \Cal H_i$.
Consider the infinite tensor product $\rho = \bigotimes_{i \in I}\rho_i$ 
of the family $(\rho_i)_{i \in I}$ 
with respect to the family $(\eta_i)_{i \in I}$. 
Recall that  $\rho$ is the representation of $S$ 
defined on the  infinite tensor product 
$\Cal H = \otimes_{i \in I}(\Cal H_i,\eta_i)$ 
of the  family of Hilbert spaces $(\Cal H_i)_{i \in I}$ 
with respect to  the family $(\eta_i)_{i \in I}$ by 
$$
\rho \Big( (\gamma_i)_{i \in I} \Big) 
\Big(
\big( \bigotimes_{f \in F}\xi_f \big)
\otimes 
\big( \bigotimes_{i \in I \setminus F} \eta_i \big)
\Big)
\, = \,
\big( \bigotimes_{f \in F} \rho_i(\gamma_i)\xi_f \big)
\otimes
\big( \bigotimes_{i \in I \setminus F} \eta_i \big),
$$
for every finite subset $F$ of $I$,
element $(\gamma_i)_{i \in I} \in S$ 
with $\gamma_i = 1$ whenever $i \in I \setminus F$,  
and decomposable vector
$(\xi_f)_{f \in F} \in \bigotimes_{f \in F} \Cal H_f$.
The representation $\rho$ is irreducible, since the
$\rho_i$'s are irreducible. 
For all this, see for example 
\cite{Guich--66}, in particular Corollary~2.1.
\par

Let us check that, 
for $\gamma = \left( \gamma_i \right)_{i \in I} \in S$, $\gamma \ne e$,
the operator $\rho(\gamma)$ is not a multiple of the identity operator.
Choose $j \in I$ such that $\gamma_j \ne e$.
Observe that the set
$$
\{\delta \in S_j \ : \ 
\text{$\rho_j(\delta)$ is  a multiple of the identity operator} \}
$$
is an abelian normal subgroup of $S_j$ 
and is therefore reduced to $\{e\}$
since $S_j$ is simple and nonabelian. 
The operator $\rho_j (\gamma_j)$ is therefore 
not a multiple of the identity operator. 
We can now write 
$$
\Cal H \, = \, \Cal H_{j} \otimes \Cal H'_{j} 
\qquad \text{and} \qquad
\rho \, = \,  \rho_{j} \otimes \rho'_{j}
$$
where $\rho'_{j}$ is the tensor product of 
the family $(\rho_{\ell})_{\ell \in I \setminus \{j\}}$,
defined
on $\Cal H'_{j} =  \bigotimes_{\ell \in I \setminus \{j\}}
(\Cal H_{\ell},\eta_{\ell})$. 
Lemma~12 implies that  $\rho(\gamma)$  
is not a multiple of the identity operator. 
\hfill $\square$ 
\enddemo

It remains to show that $(ii) \Longleftrightarrow  (iv)$
in Theorem 2.
This will be a consequence of the following lemma.
\par
We are most grateful to Roland L\"otscher, 
who pointed out a mistake at this point in a first version of our paper;
we are also grateful to Jacques Th\'evenaz
for a helpful discussion on modular representations.

\proclaim{14.~Lemma} 
Let $\Gamma$ be a countable group; set $A=MA(\Gamma)$. 
Let $\{A_i\}_{i \in I}$ be a set of finite abelian feet of $\Gamma$
as in Proposition~1, so that $ A = \bigoplus_{i \in I} A_i$.
For each finite subset $E$ of ~$I$, set $A_E = \bigoplus_{i \in I} A_i$,
which is a finite abelian group.
Let $\widehat A, \widehat A_E$ denote the dual group of  $A, A_E$ respectively. 
The following properties are equivalent:
\roster
\item"(i)"
$A$ has a $\Gamma$--faithful character;
\item"(ii)"
there exists a character 
$\chi \in \widehat A$ such that the subgroup generated by  
$\chi^{\Gamma} \Doteq 
\{\chi^\gamma \ \mid \ \gamma \in \Gamma\}$ is dense in $\widehat A$;
\item"(iii)"
for every finite subset $E$ of $I$, 
the finite group $\widehat A_E$ has a $\Gamma$--faithful character.
\item "(iv)"
for every finite subset $E$ of $I$, there exists $\chi\in \widehat A_E$
such that $\widehat A_E$ 
is generated by the $\Gamma$--orbit  of  $\chi$.
\item"(v)"
for every finite subset $E$ of $I$, 
there exists  $x_E \in A_E$ 
such that $A_E$ is generated by the $\Gamma$--conjugacy class of  $x_E$.
\endroster 
\endproclaim

\demo{Proof}
{\it Equivalence of $(i)$ and $(ii)$ and equivalence of $(iii)$ and $(iv)$.}
Let $N$ be a normal abelian subgroup of $\Gamma$. Let $\chi \in \widehat N$.  
Denote by $H$ the closed subgroup of $\widehat N$
generated by  $\chi^{\Gamma}$.
By Pontrjagin duality,
the unitary dual of the compact abelian group $\widehat N/H$
can be identified with the subgroup 
$$
H^{\perp} \, = \, \{ a \in N \ :\ \psi(a) = 1  
\qquad \text{for all} \qquad
\psi \in H \};
$$
observe that 
$$
H^{\perp} \, = \, \{ a\in N \ : \ \psi(a) \, = \, 1
\quad \text{for all} \quad 
\psi \in \chi^{\Gamma} \} \, = \, 
\bigcap_{\gamma\in \Gamma}  \ker(\chi^\gamma).
$$
Thus $\chi^{\Gamma}$ is dense in $\widehat N$
if and only if $H^{\perp} = \{e\}$,
namely if and only if $\chi$ is $\Gamma$--faithful.

\medskip

{\it Equivalence of (ii) and (iii).}
It is clear that $(ii)$ implies $(iii)$.
Let us assume that $(iii)$ holds;
we have to check that this implies $(ii)$. 
For every finite subset $E$ of $I$, 
denote by  $p_E: \widehat A \to \widehat A_E$ 
the canonical projection.
Consider the subset 
$$
X_E \, = \, \{
\chi \in \widehat A
\, \mid \,
\text{the $\Gamma$--orbit of  $p_E(\chi)$
generates $\widehat A_E$}
\}.
$$
Since the group $\widehat A_E$ is finite,
the subset $X_E$ of $\widehat A$ is closed.
For a finite family $E_1,\hdots,E_k$ of finite subsets of $I$,
the intersection $X_{E_1} \cap \cdots \cap X_{E_k}$
contains $X_{E_1 \cup \cdots \cup E_k}$.
By Condition~(iii), $X_E$ is non empty for any finite subset $E$ of $I$.
Since $\widehat A$ is compact, it follows that
$$
\bigcap_E X_E\neq \emptyset,
$$
where $E$ runs over all finite subsets of $I$.
Let $\chi \in \bigcap_E X_E$. It is easily checked that
$\chi$ is $\Gamma$--faithful.

\medskip


\medskip

{\it Equivalence of (iv) and (v).}
Consider a finite subset $E$ of $I$.
Recall that each $A_i$ is a finite dimensional vector space
over a prime  field $\bold F_{p_i}$, for a prime number $p_i$.
For each prime $p$, denote by $V_p$ the direct sum of
those $A_i$ with $i \in E$ which are vector spaces over $\bold F_p$,
and denote by $P$ the set of primes $p$ such that $V_p \ne \{0\}$.
We have $A_E = \bigoplus_{p \in P} V_p$.
Since the $V_p$'s are subgroups of $A_E$ of pairwise coprime orders,
every subgroup $H$ of $A_E$ is a direct sum 
$\bigoplus_{p \in P} (H \cap V_p)$.
The dual group $\widehat V_p$ of $V_p$
can be identified with the dual vector space $V_p^*$;
as before, each subgroup $H^*$ of $V_p^*$
is a direct sum $\bigoplus_{p \in P} (H^* \cap V_p^*)$.
It follows that, in order to prove the equivalence
of  $(iv)$ and $(v)$, we can assume 
that $P$ consists of a single element $p$.
We can also assume that $\Gamma$ is a subgroup
of $GL(V_p)$.

Let $\bold F_{p}[\Gamma]$ denote 
the group algebra of $\Gamma$ over $\bold F_{p}$.
Observe that $V_p$ is  a semi--simple $\bold F_{p}[\Gamma]$--module,
since $V_p$ is a direct sum of minimal normal subgroups of $\Gamma$.
(A module is {\it semi--simple} if it is a direct sum of simple modules;
other authors use the terminology {\it completely reducible}.)

Under the identification of $\widehat V_p$ with $V_p^*$, 
the  $\Gamma$--action on $\widehat V_p$ corresponds to 
the dual (or contragredient)  action of $\Gamma$ on $V_p^*$.
Observe  that $V_p^*$ is  a semi--simple $\bold F_{p}[\Gamma]$--module.
Indeed, if $W$ is  submodule of  $V_p^*$, 
then its annihilator $W^{\perp}$ has a complement $Z$ in 
$V_p$ and  $Z^\perp$ is a complement of $W$ in $V_p^*$
(compare with Lemma 6.2 in \cite{Landr--83}).

Observe also that there exists $x\in V_p$
such that $V_p$ is generated by the $\Gamma$--conjugacy class  of  $x$
(respectively, there exists $\chi\in \widehat V_p$
such that $\widehat V_p$ 
is generated by the $\Gamma$--orbit  of  $\chi$)
if and only if  $V_p$ (respectively $V_p^*$)
is isomorphic, as $\bold F_{p}[\Gamma]$--module,
to a quotient of the left regular module  $\bold F_{p}[\Gamma]$.
To conclude the proof, we show that $V_p$
is isomorphic to  a quotient of $\bold F_{p}[\Gamma]$
if and only if $V_p^*$
is isomorphic to a quotient of $\bold F_{p}[\Gamma]$.

We  first show that every semi--simple submodule  of 
$\bold F_{p}[\Gamma]$ is isomorphic
to a quotient of $\bold F_{p}[\Gamma]$.
Indeed, let $\bold F_{p}[\Gamma]=\bigoplus_{j\in J} P_j$
be a direct sum decomposition of $\bold F_{p}[\Gamma]$
into  indecomposable submodules $P_j$.
Every $P_j$ contains a unique simple module $S_j$. 
Moreover, $S_j$ is isomorphic to a quotient 
 of $P_j$ and $M=\bigoplus_{j\in J} S_j$ is the 
sum of all simple submodules of $\bold F_{p}[\Gamma]$.
For  the standard facts on representation theory
of finite groups, see for example
\cite{Landr--83}, in particular Theorem 6.8.
Let $N$ be a semi--simple submodule  of 
$\bold F_{p}[\Gamma]$. 
Then $N$ is a submodule of $M$ and  is therefore isomorphic
to $\bigoplus_{j\in J'} S_{j}$ 
for a subset $J'$ of $J$.  Hence,
$N$ is isomorphic
to a quotient of $\bigoplus_{j\in J'} P_{j}$.
Since $\bigoplus_{j\in J'} P_{j}$ is a direct summand of 
$\bold F_{p}[\Gamma]$, 
it follows that $N$ is isomorphic
to a quotient of $\bold F_{p}[\Gamma]$
and this proves our claim.

Assume that $V_p$ is isomorphic
to a quotient of $\bold F_{p}[\Gamma]$.
Then $V_p^*$ is isomorphic
to a submodule of $\bold F_{p}[\Gamma]^*$.
Now, it is standard that 
$\bold F_{p}[\Gamma]^*$ is isomorphic to $\bold F_{p}[\Gamma]$
as a $\bold F_{p}[\Gamma]$--module (see
\cite{Landr--83, Theorem 6.3}).
Hence, $V_p^*$ is isomorphic to a submodule of $\bold F_{p}[\Gamma]$.
By what we have seen above, 
it follows that $V_p^*$  is   isomorphic to a quotient of $\bold F_{p}[\Gamma]$.
Similarly, if $V_p^*$  is   isomorphic to a quotient of $\bold F_{p}[\Gamma]$,
then  $V_p$ is isomorphic to a quotient of $\bold F_{p}[\Gamma]$.
\hfill $\square$ 
\enddemo

\bigskip
\head{\bf
4.~Finite and infinite dimensional representations
}\endhead

Our proof of Theorem~4 uses the following elementary lemma,
which is well--known. 
To our surprise, 
we haven't been able to find a convenient reference.

\proclaim{15.~Lemma}
Let $\Omega$ a standard Borel space and
$\mu$ a bounded positive measure  on $\Omega$.
\par

(i) Let $A$ be a separable C$^*$--algebra, 
$\underline\pi$ a representation of $A$, and
$$
\underline\pi \, = \, \int^{\oplus}_{\Omega} \underline\pi_{\omega} d\mu(\omega)
$$
a direct integral decomposition of $\underline\pi$
with respect to a measurable field $\omega \longmapsto \underline\pi_{\omega}$
of representations of $A$.
Then $\underline\pi_{\omega}$ is weakly contained in $\underline\pi$ 
for almost all $\omega$ in $\Omega$.
\par

(ii) Let $\Gamma$ be a countable group,
$\pi$ a representation of $\Gamma$, and
$$
\pi \, = \, \int^{\oplus}_{\Omega} \pi_{\omega} d\mu (\omega)
$$
a direct integral decomposition of $\pi$
with respect to a measurable field $\omega \longmapsto \pi_{\omega}$
of representations of $\Gamma$.
Then $\pi_{\omega}$ is weakly contained in $\pi$
for almost all $\omega \in \Omega$.
\endproclaim

\demo{Proof}
(i) By definition of \lq\lq weak containment\rq\rq ,
we have to show that 
$\ker (\underline\pi) \subset \ker({\underline\pi}_{\omega})$
for almost all $\omega \in \Omega$.
Since $A$ is separable, so is $\ker(\underline\pi)$,
and we can choose in this kernel a countable dense subset, 
say $C$.
For any $x \in A$, recall from the theory of direct integrals
that $\Vert \underline\pi(x) \Vert$ 
is the essential supremum (on $\omega \in \Omega$)
of the norms $\Vert \underline\pi_{\omega}(x) \Vert$,
so that $\Vert \underline\pi_{\omega}(x) \Vert \le \Vert \underline\pi(x) \Vert$
for almost all $\omega \in \Omega$;
in particular, any $x \in C$ is in $\ker(\underline\pi_{\omega})$ 
for almost all $\omega \in \Omega$.
Since $C$ is countable, we have also
$C \subset \ker(\underline\pi_{\omega})$ 
for almost all $\omega \in \Omega$,
and this implies the announced conclusion.
\par

(ii) Any representation $\pi$ of $\Gamma$
corresponds to a representation $\underline \pi$ 
of the maximal C$^*$--algebra 
$A = C^*_{\operatorname{max}}(\Gamma)$ of the group.
For two representations $\pi_1, \pi_2$ of the group,
$\pi_1$ is {\it weakly contained} in $\pi_2$ 
if and only if 
$\ker(\underline\pi_2) \subset \ker(\underline\pi_1)$;
moreover, a direct integral decomposition
$\pi = \int^{\oplus}_{\Omega} \pi_{\omega} d\mu (\omega)$ 
at the level of $\Gamma$
corresponds to a direct integral decomposition
$\underline\pi = \int^{\oplus}_{\Omega}  \underline\pi_{\omega} d\mu (\omega)$
at the level of $C^*_{\operatorname{max}}(\Gamma)$,
with the same space $\Omega$ and the same measure~$\mu$.
Thus (ii) is a consequence of (i).
\par

[More generally, both (ii) and its proof hold {\it verbatim} for
representations of separable locally compact groups.]
\hfill $\square$
\enddemo

To prove Theorem~4, it is clearly enough to show that
Conditions (i) and (ii) there are equivalent.
The implication $(i) \Longrightarrow (ii)$
is a straightforward consequence
of \cite{Thoma--64, Korollar~1}, according to which {\it every}
irreducible representation of a virtually abelian group 
is finite dimensional.

\demo{End of proof of Theorem~4, namely of $(ii) \Longrightarrow (i)$}
We assume that $\Gamma$ has Property~(ii), 
and we split the proof in two cases.
\par

Assume first that $\Gamma$ is not amenable. 
Let $\sigma$ be a $\Gamma$--faithful irreducible representation of $MS(\Gamma)$;
set $\pi = \operatorname{Ind}_{MS(\Gamma)}^{\Gamma} \sigma$.
By Lemmas~10 and~15, some (in fact almost every) 
irreducible representation $\pi_0$
which occurs in some direct integral decomposition of $\pi$
is faithful and is weakly contained in $\pi$.
As $MS(\Gamma)$ is amenable, $\pi$ is weakly contained
in the left regular representation of~$\Gamma$,
and therefore the same holds for $\pi_0$.
As $\Gamma$ is not amenable, $\pi_0$ cannot be finite dimensional,
so that $\Gamma$ has Property~(i).
\par

Assume now that $\Gamma$ is amenable.
Assume furthermore, by contradiction, 
that $\Gamma$ does not have Property~(i).
Then $\Gamma$ has a finite dimensional faithful irreducible representation,
by the first part of~(ii).
In particular, $\Gamma$ can be viewed
as a subgroup of the compact unitary group $\Cal U(n)$,
for some integer $n \ge 1$.
By Tits' alternative \cite{Tits--79}, there exists in $\Gamma$
a soluble subgroup $\Delta$ of finite index.
Let $R$ denote the closure of $\Delta$ in $\Cal U(n)$
and let $R^0$ denote its connected component;
then $R^0$ is of finite index in $R$ 
(because $R$ is a compact Lie group, see for example 
\cite{Helga--62, Chapter II, Theorem 2.3})
and an abelian group (because a connected compact  group is soluble
if and only if it is abelian, see for example \cite{Bourb--82, Appendice~I}).
Thus $\Delta \cap R^0$ is an abelian subgroup of finite index in $\Gamma$;
but this contradicts the hypothesis that $\Gamma$ is not virtually abelian,
and this ends the proof.  \hfill $\square$
\enddemo

\proclaim{16.~Proposition}
Let $\Gamma$ be a countable group.
\par

(i) If there exists a factor $M$ and an injective homomorphism
$\pi : \Gamma \longrightarrow \Cal U (M)$ such that $\pi(\Gamma)'' = M$,
then $\Gamma$ is irreducibly represented.
\par

(ii) If $\Gamma$ is irreducibly represented, 
then there exists a factor $M = \Cal L (\Cal H)$  of type I
and a faithful representation 
$\pi : \Gamma \longrightarrow \Cal U (\Cal H)$
such that $\pi(\Gamma)'' = \Cal L (\Cal H)$.
\endproclaim

\demo{Proof}
Let $\pi$ be as in (i). 
If $M$ is an algebra of operators on some Hilbert space $\Cal K$,
then $\pi$ is in particular 
a factorial representation of $\Gamma$ in $\Cal K$.
It corresponds to a C$^*$--representation, 
say $\underline{\pi} : C^*_{\operatorname{max}}(\Gamma)
\longrightarrow \Cal L (\Cal K)$.
By a result of Dixmier (Corollary 3, Page 100 of \cite{Dix--60}),
there exists an irreducible representation $\underline{\rho}$
of $C^*_{\operatorname{max}}(\Gamma)$ 
such that $\underline{\pi}$ and $\underline{\rho}$
have the same kernel.
The restriction $\rho$ of $\underline{\rho}$ to $\Gamma$ 
is therefore a faithful irreducible representation.
\par

In view of Schur's lemma, (ii) is nothing but a reformulation of the definition
of \lq\lq irreducibly represented\rq\rq .
\hfill $\square$
\enddemo

Corollary~5 is a straightforward consequence 
of Theorem~4 and Proposition~16.
\par

Short of knowing how to answer the questions which follow Corollary~5, 
let us record the following elementary remark.

\proclaim{17.~Observation}
If $\Gamma$ is a  countable group 
which has a finite dimensional faithful irreducible representation,
then $MH(\Gamma)$ is a finite group.
\endproclaim

\demo{Proof}
Consider the following properties
of a group $\Gamma$:
\roster
\item"(a)" 
$\Gamma$ has a finite dimensional faithful irreducible representation;
\item"(b)" 
$MS(\Gamma)$ has a finite dimensional $\Gamma$--faithful irreducible
representation;
\item"(c)"
$MA(\Gamma)$ has a $\Gamma$--faithful character and
$MH(\Gamma)$ has a finite dimensional faithful irreducible representation;
\item"(d)"
$MA(\Gamma)$ has a $\Gamma$--faithful character and
$MH(\Gamma)$ is a finite group.
\endroster 
Property~(a) implies Property~(b) by Lemma~9,
Properties (b) and (c) are equivalent because
$MS(\Gamma) = MA(\Gamma) \oplus MH(\Gamma)$,
and Properties (c) and (d) are equivalent because
$MH(\Gamma)$ is a direct sum of finite simple groups.
\par
[Observe that, however, Property~(b) does not imply Property~(a):
if $\Gamma$ is an icc group which does not have any finite dimensional
faithful representation, 
for example the group of permutations of finite support
of $\bold Z$, 
then $\Gamma$ has Property (b) since $MS(\Gamma) = \{e\}$,
but does not have Property~(a).]
\hfill $\square$
\enddemo

About Conditions (ii) and (iv) of Theorem~4,
let us moreover recall the following facts.
For countable groups, 
and more generally for separable locally compact groups
and for separable C$^*$--algebras,
there is a notion of being {\it of type~I},
defined in terms of the von Neumann algebras
generated by the images of appropriate representations.
It is then a theorem of Thoma that a countable group is of type~I
if and only if it is virtually abelian,
if and only if
all its irreducible  representations are finite dimensional.
See \cite{Thoma--64} and \cite{Glimm--61}.


\bigskip
\head{\bf
5.~Remarks
}\endhead

\subhead 
5.a.~Minisocles, socles, and examples 
\endsubhead

The {\it socle} of a group $\Gamma$ is the subgroup $S(\Gamma)$
generated by the union of the minimal normal subgroups (finite or infinite).
Here are some examples of socles and minisocles.
\smallskip

(I) For a prime $p$ and an integer $n \ge 1$, 
the socle of the finite cyclic group $\bold Z /p^n \bold Z$
is isomorphic to $\bold Z /p\bold Z$.
The socle of the finite symmetric group $\operatorname{Sym}(n)$
is the corresponding alternating group if $n = 3$ or $n \ge 5$,
and the Vierergruppe if $n = 4$.
\par

If $\Gamma$ is a $2$--transitive subgroup of $\operatorname{Sym}(n)$,
then $S(\Gamma)$ is either of the form $(\bold F_p)^m$ 
or a finite simple group. More generally and more precisely, 
if $\Gamma$ is a primitive subgroup of $\operatorname{Sym}(n)$,
the O'Nan--Scott Theorem (1980) provides detailed informations
on the socle of~$\Gamma$; in particular, $S(\Gamma) \approx S^m$
for some finite simple group $S$ and some integer $m$.
See for example Chapter 4 in \cite{DixMo--96}.

\smallskip

(II) Free abelian groups $\bold Z^n$, $n \ge 1$, 
and nonabelian free groups have socles reduced to one element.
\smallskip

(III) For $n \ge 3$, the socle of $SL_n(\bold Z)$ 
is reduced to one element or of order two, 
if $n$ is odd or even respectively
(because any noncentral normal subgroup of $SL_n(\bold Z)$ 
contains a congruence subgroup, and consequently is never minimal).
\smallskip

(IV) Let $\Gamma$ be a lattice in a connected semisimple Lie group $G$
with finite center $Z(G)$ and without compact factor.
It is an easy consequence of the Borel density theorem
that, if the centre of $\Gamma$ is $\{e\}$, then $\Gamma$ is icc,
so that $MS(\Gamma) = \{e\}$;
more generally, $MS(\Gamma) = \Gamma \cap Z(G)$.
\par

The minisocle of a just infinite group 
is reduced to one element (by definition). 
In particular, the minisocle of the Grigorchuk group is reduced to one element.
\smallskip

(V) If $\Gamma$ is a direct sum 
of a  family of infinite simple groups, 
then $MS(\Gamma) = \{e\}$ and $S(\Gamma) = \Gamma$.
If $\Gamma$ is a direct sum 
of a  family of finite simple groups, 
then $MS(\Gamma) =  \Gamma$.
\smallskip

(VI) The socle of a nilpotent group $\Gamma$ 
is contained in the centre $Z(\Gamma)$ of $\Gamma$,
because $N \cap Z(\Gamma) \ne \{e\}$
for any normal subgroup $N \ne \{e\}$ of $\Gamma$.
\smallskip

(VII) Let $\Gamma$ be an abelian torsion-free group
with cardinal strictly larger than that of the real numbers,
for example a direct product of copies of $\bold Z$
indexed by $\bold R$.
Then $MS(\Gamma)$ is reduced to one element,
and $\Gamma$ does not have any faithful character,
so that the equivalences of Theorem~2 do not hold for $\Gamma$.
\smallskip

(VIII) Let $H$ be a group, $p$ a prime number,
$U$ a vector space over the prime field with $p$ elements,
$\pi : H \longrightarrow GL(U)$ a faithful representation
which is semi--simple (namely which is a direct sum of irreducible representations),
and $\Gamma = H \ltimes U$ the corresponding  semi--direct product.
Then $U$ is the socle of $\Gamma$.
\par
Indeed, let $N$ a minimal normal subgroup of $\Gamma$.
If $N \cap U \ne \{0\}$, then $N \subset U$,
and moreover $N$ is a $H$--invariant subspace of $U$
wich is irreducible,  by minimality; these $N$'s generate $U$.
If one had $N \cap U = \{0\}$, then $N$ and $U$ would commute
(being two normal subgroups of $\Gamma$), so that $N$
would act trivially on $U$, and this is ruled out 
by the faithfulness of $\pi$.
\par
Let $U$ be of the form 
$U = \left( \bigoplus_{i \in I}V_i \right) \oplus \left( \bigoplus_{j \in J} W_j \right)$,
with each $V_i$ a $H$--invariant irreducible finite--dimensional subspace of $U$, 
and each $W_j$ a $H$--invariant irreducible infinite--dimen\-sional subspace of $U$.
Then the mini--socle of $\Gamma$ is $\bigoplus_{i \in I} V_i$.
\par
The construction carries over to the case where each $V_i$ and $W_j$
is a vector space over a prime field of which the number of elements
depends on $i$ and $j$.

\medskip
\subhead
5.b.~Minisocles, FC--kernels, and P.~Hall's theorems
\endsubhead

The {\it FC--kernel} of a group $\Gamma$ 
is the subset $\Gamma_{\operatorname{FC}}$ 
of $\Gamma$  of elements which have a finite conjugacy class. 
It is a characteristic subgroup of $\Gamma$.
\par

The {\it periodic FC--kernel} of $\Gamma$ 
is the subset $\Gamma_{\operatorname{FC}}^{\operatorname{per}}$ 
of $\Gamma_{\operatorname{FC}}$ of elements of finite order.
It is also a subgroup of $\Gamma$,
indeed a locally finite subgroup (Dicman's Lemma, see e.g. \cite{Tomki--84}).
It follows from the definitions and from Dicman's Lemma that
$$
MS(\Gamma)
\, \subset \,
\Gamma_{\operatorname{FC}}^{\operatorname{per}}
$$
(the inclusion can be strict, as it is for example the case
if $\Gamma$ is cyclic of order four).
\par

Any subgroup of a restricted direct product of finite groups 
is a periodic FC--group which is residually finite, 
and any quotient of a periodic FC--group is a periodic FC--group.
For countable groups, Philip Hall has established in 1959 
the converse implications:
\roster
\item""
any {\it countable} periodic FC--group which is residually finite 
can be embedded in a restricted direct product of finite groups,
and any {\it countable} periodic FC--group is isomorphic to 
a quotient of a subgroup of a restricted direct product of finite groups
\endroster
(Theorems 2.5 and 3.2 in \cite{Tomki--84}).

\medskip

A {\it hoof} of a group $\Gamma$ is a foot of a foot.
Thus, with the notation of Proposition~1, the subgroups $\bold F_p$
and $S_1$  are hooves of $\Gamma$.
\par

Let $\Gamma$ be a group which has a finite Jordan--H\"older sequence
(for example a finite group);
if a simple group $S$ is a foot of $\Gamma$,
then $S$ is isomorphic to a quotient of some Jordan-H\"older sequence of $\Gamma$
(Bourbaki, Alg\`ebre, nouvelle \'edition, 1970, chap.~I, \S~4, no~7).
But the converse does not hold: the group of order~$3$ 
is a simple quotient of a Jordan--H\"older sequence of the alternating group
$\operatorname{Alt}(4)$ of order $12$, 
but $\operatorname{Alt}(4)$ has a unique foot 
which is the Vierergruppe, of order $4$.

\medskip
\subhead
5.c. Recall of a theorem of Gelfand and Raikov
\endsubhead

Recall the following basic result of the theory of group representations,
due to Gelfand and Raikov 
(see \cite{GelRa--42}, as well as Corollary 13.6.6 in \cite{Dix--69C$^*$}):
\roster\item""
{\it
for any $\gamma \in \Gamma$, $\gamma \ne e$,
there exists an irreducible representation $\pi_{\gamma}$
such that $\pi_{\gamma}(\gamma) \ne \pi_{\gamma}(e)$.
}
\endroster
This holds for any group $\Gamma$, countable or not;
indeed, this holds for any locally compact group, 
with $\pi_{\gamma}$ a continuous representation.
There are two main ingredients of the proof:
the group has a faithful representation 
which is the left--regular representation,
and any representation has some description in terms
of irreducible representations
(via  functions of positive types and
a theorem of the Krein--Milman type).
\par

For a countable group $\Gamma$
which has the properties (ii) to (v) of Theorem~2,
we have shown that $\pi_{\gamma}$ can be chosen 
{\it independently of $\gamma$.}

\medskip
\subhead
5.d. Recall on tensor powers of faithful representations
\endsubhead

Let $\Gamma$ be a group
and let $\pi$ be a faithful representation of $\Gamma$.   
For integers $m,n\geq 0,$ consider the tensor power 
$\pi_{m,n}=\pi^{\otimes m} \otimes \overline{\pi}^ {\otimes n}$,
where $ \overline{\pi}$ denotes the  representation conjugate to $\pi$
and $\pi^{\otimes m}$ the  tensor product of $m$ copies of $\pi$.
Then the left regular representation of $\Gamma$ 
is weakly contained in the direct sum  $\bigoplus_{m,n \geq 0}\pi_{m,n}$
(see Example 1.11 in [BeLaS--92]).
This is a generalisation,
with weak containment replacing strong containment,
of a well--known fact about finite groups 
(and compact groups, 
see \cite{Cheva--46}, Chapter~VI, \S~VII, Proposition~3).
Thus, if, in addition, $\Gamma$ is amenable, 
then every  representation of $\Gamma$
is weakly contained in $\bigoplus_{m,n \geq 0}\pi_{m,n}$. 
(All this carries over to locally compact groups.)
\par

For a countable group $\Gamma$
which has the properties (ii) to (v) of Theorem~2,
the representation $\pi$ can be chosen to be
faithful {\it and irreducible.}

\medskip
\subhead
5.e.~Primitive group C$^*$--algebras
\endsubhead

Denote by $C^*_{\operatorname{red}}(\Gamma)$ the reduced C$^*$--algebra,
and as above by $C^*_{\operatorname{max}}(\Gamma)$ the maximal C$^*$--algebra
of a group $\Gamma$. A representation of either one of these algebras
is irreducible if and only if its restriction to $\Gamma$ is irreducible.
It follows that, 
{\it if one of 
$C^*_{\operatorname{red}}(\Gamma)$, $C^*_{\operatorname{max}}(\Gamma)$
is primitive, then $\Gamma$ is irreducibly represented.}
\par

Many examples of countable groups are known for which
$C^*_{\operatorname{red}}(\Gamma)$ is simple \cite{Harpe--07}, 
and  {\it a fortiori} primitive.
These groups are in particular irreducibly represented.
Concerning the properties of $\Gamma$ and $C^*_{\operatorname{red}}(\Gamma)$, 
consider the three following conditions:
\roster
\item"(NF)" 
   $\Gamma$ does not have any finite normal subgroup besides $\{e\}$;
\item"(NA)" 
   $\Gamma$ does not have any amenable normal subgroup besides $\{e\}$;
\item"(C$^*$S)" $C^*_{\operatorname{red}}(\Gamma)$ is simple.
\endroster
It is straightforward that (NF) is a rephrasing of 
the condition $MS(\Gamma) = \{e\}$, 
and that it follows from (NA).
It is elementary to check that (NA) follows from (C$^*$S),
but we recall that it is not known whether the converse holds
(see \cite{BekHa--00} and \cite{Harpe--07}).
\par

If $\Gamma$ is amenable, the C$^*$--algebras
$C^*_{\operatorname{red}}(\Gamma)$ and $C^*_{\operatorname{max}}(\Gamma)$
are isomorphic.
They are primitive if and only if $\Gamma$ is icc \cite{Murph--03}.
\par

If $\Gamma$ is a nonabelian free group,
it is a  result of Yoshiwaza
that $C^*_{\operatorname{max}}(\Gamma)$ is primitive
(see  \cite{Yoshi--51}, as well as \cite{Choi--80}).
See the discussion around Problem~25 in \cite{Harpe--07}.

\bigskip
\head{\bf
6.~A generalisation of Theorem~2 
}\endhead
Consider a countable group $\Gamma$ 
and a subgroup $G$ of the automorphism group of $\Gamma$
which contains all inner automorphisms.
There is an obvious notion of $G$--faithful representation,
which coincides with that of $\Gamma$--faithful representation
in case $G$ coincides with the group of inner automorphisms.
Observe that $MS(\Gamma)$, $MA(\Gamma)$, and $MH(\Gamma)$
are $G$--invariant subgroups of $\Gamma$,
since all three are characteristic.
\par

We define a {\it $G$--foot} to be 
a minimal $G$--invariant subgroup of $\Gamma$.
Let $\Cal F^G_{\Gamma}$ denote the set of finite $G$--feet of $\Gamma$;
it is the union of the set $\Cal A^G_{\Gamma}$ of abelian finite $G$--feet
and of its complement $\Cal H^G_{\Gamma}$.
The {\it $G$--minisocle} of $\Gamma$ is the subgroup $MS^G(\Gamma)$
generated by its $G$--feet,
and we have as in Section~1 subgroups 
$MA^G(\Gamma)$ and $MH^G(\Gamma)$.

\proclaim{18.~Proposition} Let $\Gamma$ and $G$ be as above.
\par

(i) Each $B \in \Cal A^G_{\Gamma}$ is isomorphic to $(\bold F_p)^n$
for some prime $p$ and some positive integer $n$ (depending on $B$).
\par

(ii) There exists a subset $\{B_i\}_{i \in I}$ of $\Cal A^G_{\Gamma}$
such that $MA^G(\Gamma) = \bigoplus_{i \in I} B_i$.
In particular, the group $MA^G(\Gamma)$ is abelian.
\par

(iii) For each $H \in \Cal H^G_{\Gamma}$,
the feet $S_1,\hdots,S_k$ of $H$ 
are conjugate under $G$, and simple. 
Moreover $H = S_1 \oplus \cdots \oplus S_k$.
\par

(iv) We have $MH^G(\Gamma) = \bigoplus_{H \in \Cal H^G_{\Gamma}} H$.
\par

(v) We have $MS^G(\Gamma) = MA^G(\Gamma) \oplus MH^G(\Gamma)$.
\endproclaim

\proclaim{19.~Theorem} Let $\Gamma$, $G$, and
 $MA^G(\Gamma) = \bigoplus_{i \in I} B_i$ be as above.
The following properties are equivalent:
\roster
\item"(i)" 
$\Gamma$ has a representation which is irreducible and $G$--faithful;
\item"(ii)" 
$MA(\Gamma)$ has a $G$--faithful character;
\item"(ii')"
$MA^G(\Gamma)$ has a $G$--faithful character;
\item"(iii)" 
$MS(\Gamma)$ has a $G$--faithful irreducible representation;
\item"(iii')"
$MS^G(\Gamma)$ has a $G$--faithful irreducible representation;
\item"(iv)" 
for every finite subset $E$ of $I$, 
there exists an element
$x_E$ in $MA^G_E(\Gamma) \Doteq \bigoplus_{i \in E}B_i$  such that
the $G$--orbit of $x_E$ generates $MA^G_E(\Gamma)$.
\endroster
In particular, a countable group $\Gamma$ 
has a $G$--faithful irreducible representation
as soon as $MA^G(\Gamma) = \{e\}$, 
{\it a fortiori} as soon as $MS^G(\Gamma) = \{e\}$.
\endproclaim

For example, let $\Gamma = \bigoplus_{i \in \bold N} A_i$ be a countable infinite
direct sum of groups $A_i$ indexed by the natural numbers, 
each of them isomorphic to a given finite cyclic group,
and let $G$ be the group of permutations of $\bold N$,
identified in the natural way to a group of automorphisms of $\Gamma$.
Then $\Gamma$ is irreducibly underrepresented, but has a $G$--faithful
irreducible character, for example the projection onto $A_1$ followed by the
natural isomorphisms of $A_1$ with the appropriate group of roots of unity.

\bigskip

Proposition 18 and Theorem 19 can be proved by essentially the same
arguments as in Sections 2 and 3. 
\par

Lemmas 9 and 10 should be reformulated 
for a $G$--invariant subgroup $N$ of $\Gamma$; 
in the new Lemma 9, the $G$--faithfulness of $\pi$ implies that 
$\sigma_{\omega}$ is $G$--faithful
for almost all $\omega \in \Omega$; 
in the new Lemma 10, if we assume that $\sigma$ is $G$--faithful
and that there does not exist 
any finite $G$--foot $M$ such that $M \cap N = \{e\}$,
then $\pi_{\omega}$ is $G$--faithful
for almost all $\omega \in \Omega$.
In the new Lemma~13, the groups $N,A,S$ should be  $G$--invariant,
but $S_i$ normal (not necessarily $G$--invariant) and simple (not necessarily
$G$--simple); the conclusion is that $N$ has a $G$--faithul irreducible
representation if and only if $A$ has a $G$--faithful character.
In the new Lemma~14, both $A$ and the $A_i$ should be 
$G$--invariant, and $\chi^{\Gamma}$ should be replaced by~$\chi^G$.
The other (minor) modifications are left to the reader.


\bigskip
\Refs
\widestnumber\no{ReVaW--02} 

\ref \no BekHa--00 \by B. Bekka and P. de la Harpe
\paper Groups with simple reduced C$^*$--algebras
\jour Expositiones Math. \vol 18 \yr 2000 \pages 215--230
\endref

\ref \no BeLaS--92
\by B. Bekka, A. Lau, and G. Schlichting
\paper On invariant subalgebras  of the Fourier-Stiletjes algebra
of a locally compact groups
\jour Math. Ann.
\vol 294 \yr 1992 \pages 513--522
\endref

\ref \no Bourb--82 \by N. Bourbaki
\book Groupes et alg\`ebres de Lie, chapitre~9
\publ Masson \yr 1982
\endref

\ref \no Burns--11 \by W. Burnside
\book The theory of groups of finite order (2nd Edition)
\publ Cambridge University Press \yr 1911  
\endref

\ref \no Cheva--46 \by C. Chevalley
\book Theory of Lie groups I
\publ Princeton Univ. Press \yr 1946
\endref

\ref \no Choi--80 \by M.-D. Choi
\paper The full C$^*$--algebra of the free group of two generators
\jour Pacific J. Math. \vol 87 \yr 1980 \pages 41--48
\endref

\ref \no Dix--60 \by J. Dixmier 
\paper Sur les C$^*$--alg\`ebres
\jour Bulletin de la Soc. Math. France
\vol 88 \yr 1960 \pages 95--112
\endref

\ref \no Dix--69C$^*$ \by J. Dixmier 
\book Les C$^*$--alg\`ebres et leurs repr\'esentations 
\publ Gauthier--Villars \yr 1969 
\endref

\ref \no Dix--69vN \by J. Dixmier 
\book Les alg\`ebres d'op\'erateurs dans l'espace hilbertien
(alg\`ebres de von Neumann)
\publ deuxi\`eme \'edition, Gauthier--Villars \yr 1969 
\endref

\ref \no DixMo--96 \by J.D. Dixon and B. Mortimer
\book Permutation groups
\publ Springer \yr 1996
\endref

\ref \no Dixon--71 \by J.D. Dixon
\book The structure of linear groups
\publ Van Nostrand \yr 1971
\endref

\ref \no Gasch--54 \by W. Gasch\"utz
\paper Endliche Gruppen mit treuen absolut--irreduziblen Darstellungen
\jour Math. Nach. \vol 12 \yr  1954 \pages 253--255
\endref


\ref \no GelGl \by T. Gelander and Y. Glasner
\paper Countable primitive groups
\jour GAFA, Geometric and Functional Analysis \yr to appear 
\endref

\ref \no GelRa--42 \by I.M. Gelfand and D.A. Rajkov
\paper Irreducible unitary representations of locally bicompact groups
\jour Mat. Sb. \vol 13(55) \yr 1942 \pages 301--316
[= I.M. Gelfand, Collected papers, Volume II, pages 3--17]
\endref

\ref \no GlaMo--07 \by Y. Glasner and N. Monod
\paper Amenable actions, free products, and a fixed point property
\jour Bull. London Math. Soc. 
\vol 39 \yr 2007 \pages 138--150
\endref

\ref \no Glimm--61 \by J. Glimm  
\paper Type I C$^*$--algebras
\jour Annals of Math. \vol 73 \yr 1961 \pages 572--612 
\endref

\ref \no Guich--66 \by A. Guichardet
\paper Produits tensoriels infinis et 
repr\'esentations des relations d'anticommutation
\jour Ann. Sci. Ec. Normale Sup. \vol 83 \yr 1966 \pages 1--52
\endref

\ref \no Harpe--07 \by P. de la Harpe
\paper On simplicity of reduced C$^*$--algebras of groups
\jour Bull. London Math. Soc. 
\vol 39 \yr 2007 \pages 1--26
\endref

\ref \no Helga--62 \by S. Helgason
\book Differential geometry and symmetric spaces
\publ Academic Press \yr 1962
\endref


\ref \no Huppe--98 \by B. Huppert
\book Character theory of finite groups
\publ Walter de Gruyter \yr 1998
\endref

\ref \no Landr--83 \by P. Landrock
\book Finite group algebras and their modules
\publ London Math. Soc. Lecture Notes Series 84, Cambridge Univ. Press
\yr 1983
\endref


\ref \no Murph--03 \by G.J. Murphy
\paper Primitivity conditions for full group C$^*$-algebras
\jour Bull. London Math. Soc. \vol 35 \yr 2003 \pages 697--305
\endref

\ref \no P\'alfy--79 \by P.P. P\'alfy
\paper On faithful irreducible representations of finite groups
\jour Studia Sci. Math. Hungar. \vol 14 \yr 1979 \pages 95--98
\endref

\ref \no Remak--30 \by R. Remak
\paper\"Uber minimale invariante Untergruppen in der Theorie der endlichen
Gruppen
\jour  J. reine angew. Math. \vol 162 \yr 1930 \pages 1--16
\endref


\ref \no RO--IV \by F.J. Murray and J. von Neumann 
\paper On rings of operators, IV 
\jour Annals of Math. \vol 44 \yr 1943 \pages 716--808
\endref

\ref \no Thoma--64 \by E.~Thoma
\paper \"Uber unit\"are Darstellungen abz\"albarer, diskreter Gruppen
\jour Math. Ann. \vol 153 \yr 1964 \pages 111--138
\endref


\ref \no Tits--79 \by J.\ Tits 
\paper Free Subgroups in Linear Groups  
\jour J.~of Algebra \vol 20 \yr 1979 \pages 250-270
\endref

\ref \no Tomki--84 \by M.J. Tomkinson
\book FC--groups
\publ Pitman \yr 1984
\endref

\ref \no Tushe--93 \by A.V. Tushev
\paper On exact irreducible representations of locally normal groups
\jour Ukrainian Math. J. \vol 45:12 \yr 1993 \pages 1900--1906
\endref

\ref \no Yoshi--51 \by H. Yoshizawa 
\paper Some remarks on unitary representations of the free group 
\jour Osaka Math. J. \vol 3 \yr 1951 \pages 55--63 
\endref

\endRefs
\enddocument